\documentclass[12pt]{amsart}
\usepackage{amsmath,amscd}
\usepackage{geometry,amsfonts,amssymb,amsthm,txfonts,pxfonts,amscd} 
\usepackage{algorithmicx,algpseudocode}
\usepackage{algorithm}
\usepackage{refcheck,graphicx,url}
\norefnames
\nocitenames
\def\struckint{\mathop{%
\def\mathpalette##1##2{\mathchoice{##1\displaystyle##2}%
  {##1\textstyle##2}{##1\scriptstyle##2}{##1\scriptscriptstyle##2}}%
\mathpalette
{\vbox\bgroup\baselineskip0pt\lineskiplimit-1000pt\lineskip-1000pt
\halign\bgroup\hfill$}
{##$\hfill\cr{\intop}\cr\diagup\cr\egroup\egroup}%
}\limits}
\newtheorem{theorem}{Theorem}[section]

\theoremstyle{remark}

\newtheorem{question}[theorem]{Question}

\newcommand{\cx}{\mathbb{C}}
\newcommand{\integers}{\mathbb{Z}}

\newcommand{\hyps}{\mathbb{H}}
\newcommand{\reals}{\mathbb{R}}

\DeclareMathOperator{\tr}{tr}
\DeclareMathOperator{\acosh}{arccosh}

\DeclareMathOperator{\Sp}{Sp}
\DeclareMathOperator{\SL}{SL}

\DeclareMathOperator{\GL}{GL}

\DeclareMathOperator{\Isom}{Isom}
\DeclareMathOperator{\PickMatrix}{PickMatrix}
\DeclareMathOperator{\random}{random}
\DeclareMathOperator{\GenRandom}{GenRandom}
\DeclareMathOperator{\round}{round}
\DeclareMathOperator{\PickLatticeVector}{PickLatticeVector}
\DeclareMathOperator{\PickHyperbolic}{PickHyperbolic}
\DeclareMathOperator{\PickHalfplane}{PickHalfplane}
\DeclareMathOperator{\Reduce2}{Reduce2}
\DeclareMathOperator{\PSD}{PSD}

\begin{document}
\title[How to pick a matrix?]{How to pick a random integer matrix?\\ (and other questions)}

\author{Igor Rivin}
\address{Department of Mathematics, Temple University, Philadelphia}
\curraddr{Mathematics Department and ICERM, Brown University}
\email{rivin@temple.edu}
\date{\today}
\begin{abstract}
We discuss the question of how to pick a matrix uniformly (in an appropriate sense) at random from groups big and small. We give algorithms in some cases, and indicate interesting problems in others.
\end{abstract}
\thanks{The author would like to thank Nick Katz , Chris Hall, Peter Sarnak, and Hee Oh for helpful conversations. He would also  like to thank ICERM and Brown University for their hospitality and generous support}
\maketitle
\tableofcontents
\section{Introduction}
In a number of papers (see, for example, \cite{rivin2008walks,gorodnik2011splitting,fuchsrivin,rivin2012generic}) results are proved about the behavior of a typical element of a lattice in a semisimple Lie Group (for example, $\SL(n, \integers),$ where ``typical'' means picked uniformly at random from from all matrices in the group with (for example) Frobenius norm bounded above by a constant $X.$ While these results are often enlightening, what is not addressed is how one might actually pick such a matrix -- in this paper I try to address this question.

Suppose you are asked to pick a matrix uniformly at random from all the matrices  $M=\begin{pmatrix} a & b\\ c & d\end{pmatrix}$ in $\SL(2, \integers)$ such that the th Frobenius norm $\|M\|,$ defined as $\|M\| = \sqrt{\tr M M^t} = \sqrt{a^2+b^2+c^2+d^2}$ is at most $X.$ The simplest method is to pick a random matrix in $M^{2\times2}(\integers)$ satisfying the norm bound, check whether the determinant is equal to $1,$ throw it away if it is not, and return it if it is. We will describe the implementation of the function PickMatrix later, but now we note that the number of matrices in $M^{2\times2}(\integers)$ satisfying the norm bound is of order $X^4.$ On the other hand, it is known that the number of elements of $\SL(2, \integers)$ satisfying the norm bound is asymptotic to $6 X^2$ (see \cite{MR924457}), which means that the expected number of attempts before we succed is of order $O(X^2),$ which is \emph{exponential} in the size of the input (which is, roughly, $\log X$).
\begin{algorithm}
\label{naivealg}
\caption{naive algorithm for picking random elements from $\SL(2, \integers)$}
\begin{algorithmic}[1]
\Require $X$ is a real number greater than or equal to $\sqrt{2}.$
\Statex
\Function{PickSLMatrix1}{$X$}
\State $M \gets \PickMatrix(X,2)$ \Comment $\PickMatrix(n, X)$ returns a uniformly distributed  $n\times n$ integer matrix with Frobenius norm at most $X.$
\While{ $\det M \neq 1$}
\State $M\gets \PickMatrix(X, 2)$
\EndWhile
\State \Return $M$
\EndFunction
\end{algorithmic}
\end{algorithm}
Below, we will describe in detail a polynomial time algorithm for choosing a matrix from $\SL(2, \integers)$ with norm bounded above by $X.$ This algorithm is transcendental, not combinatorial, which is a little surprising. It is polynomial time, and it is an approximation algorithm, in the following sense: if in default form the biggest ratio of the probabilities of selecting matrices $A$ and $B$ is $\exp(1+ \epsilon),$ we can make the ratio $\exp(1 + \epsilon/k)$ at the cost of increasing the running time of the algorithm by a factor of $k.$ The rest of the paper is organized as follows: First, we discuss the baby version of the question (how to write the function $\PickMatrix$ in the na\"ive Algorithm \ref{naivealg}. Then we will discuss $\SL(2, \integers)$ in detail, and discuss how the algorithm may be extended to other matrix groups, including $\SL(n, \integers) $ for arbitrary $n.$
Finally, we briefly discuss the situation for finite matrix groups.
\subsection{How to produce random numbers with a given density?}
\label{fixeddensity} Suppose we have a positive function $f$ defined on the interval $(0, R],$ and we want to produce random numbers whose density at $x$ is proportional to $f(x)$ (when all we are given is a source $\random$ of \emph{uniform} random numbers on $[0, 1]$. This turns out to be easier than one might have thought, and described in Algorithm \ref{randalg}
\begin{algorithm}
\label{randalg}
\caption{Generating random numbers with a given density on $[0, R]$}
\begin{algorithmic}[1]
\Require $R$ is a real number greater than $0,$ $f$ a positive function on $[0, R].$
\Statex
\Function{GenRandom}{$f$, $R$}
\State $F(t) \gets \int_0^t f d\lambda.$ \Comment $F$ is the antiderivative of $f.$
\State $x \gets F(R)\random()$
\label{thex} \Comment Generate uniform random number between $0$ and $F(R).$
\State \Return $F^{-1}(x).$ \Comment Where $F^{-1}$ is the inverse function.
\EndFunction
\end{algorithmic}
\end{algorithm}
To show that Algorithm \ref{randalg} works, we note that  the probability that $\GenRandom{R}$ is between $t$ and $t+\Delta t$ is the probability that $x$ (on line \algref{randalg}{thex}) is between $F(t)$ and $F(t + \Delta t),$ which is about $f(t) \Delta t/F(R),$ as advertised.
\section{Geometric preliminaries}
\subsection{Uniform random points in balls}
\label{randomballs}
Suppose we want to generate a uniformly random point in a ball of radius $R$ in $\mathbb{R}^n.$ This point will have a radius and a spherical coordinate, so we generate these separately. For the spherical coordinate, it is well-known that a vector whose coordinates are identical independently distributed gaussians has direction uniformly distributed on the unit sphere (a fast in practice method of generating this is the Box-Muller method \cite{box1958note}). As for the radius, we generate a random number between $0$ and $R^n,$ then take its $n$-th root, multiplied by an appropriate constant.

Suppose now that we want to generate a random uniform point from a disk in the hyperbolic plane. The angle here is even easier (a uniform random number between $0$ and $2\pi$ will do). As for the radius, we know that the area of a disk of radius $R$ in the hyperbolic plane is $2\pi (\cosh R - 1),$ so to generate our radius, we compute a random number $x$ between $0$ and $\cosh R  - 1,$ then use $\acosh(x+1)$ as the radius. We summarize this as follows:
\begin{algorithm}
\label{pickhyperbolic}
\caption{picking a point uniformly at random from a hyperbolic disk of radius $R$}
\begin{algorithmic}[1]
\Require $X$ is a positive real number
\Statex
\Function{PickHyperbolic}{$X$}
\State $x \gets (\cosh R -1) \random()$
\State $\theta \gets 2\Pi \random()$
\Return $(\acosh(x+1), \theta)$
\EndFunction
\end{algorithmic}
\end{algorithm}

\subsection{Computing a random integer matrix} How do we write our procedure $\PickMatrix?$ The first observation is that a random $n\times n$ matrix with Frobenius norm bounded by $X$ is simply an $n^2$-tuple of  integers $a_{11}, a_{12}, \dotsc, a_{nn}$ with 
\[\sum_{i=1}^n \sum_{j=1}^n a_{ij}^2 \leq X^2,\] so we are looking for a uniformly distributed integer lattice point in the ball of radius $X$ in $\reals^{n^2}.$ The simplest (combinatorial) way to pick such a point is to pick a lattice point in the cube $[-X, X]^{n^2},$ and then throw out those points with norm bigger than $X.$ This is a perfectly fine algorithm in small dimensions, but it degrades horribly in high dimensions, since the ratio of the volume of the ball to the ratio of circumscrbed cube goes to zero superexponentially as dimension goes to infinity. In particular, for $4\times 4$ matrices, we will reject around $300000$ matrices for each one accepted. Instead, the following is an efficient algorithm:
\begin{algorithm}
\label{pickmatrix}
\caption{picking a random lattice vector of  $L^2$ norm bounded by $X$ in $\reals^n.$}
\begin{algorithmic}[1]
\Require $X$ is a positive real number.
\Function{PickLatticeVector}{$n$,$X$}
\Loop
\State $x \gets \mbox{Random vector in $\reals^n$ of norm bounded above by $X+\sqrt{n}.$}$
\label{additive}
\State $v \gets\mbox{closest lattice point to $x.$}$
\If{$\|v\| \leq X$}
\State \Return $v$
\EndIf
\EndLoop
\EndFunction

\Function{PickMatrix}{$n$, $X$}
\State \Return $\PickLatticeVector(n^2,X)$
\EndFunction
\end{algorithmic}
\end{algorithm}
Note that the additive constant of $\sqrt{n}$ (the length of a diagonal of a unit cube in $\reals^n$ (and the consequent possible resampling) is added to eliminate ``edge effect'' -- without it, the probabilities of choosing numbers close to the norm bound would be different from that of choosing smaller numbers.

\section{Action of $\SL(2, \reals)$ and $\SL(2, \integers)$ on the upper half plane}
\label{sl2action}
Recall that $\SL(2, \reals)$ acts on the upper halfplane $H=\{z \left| \Im z  > 0\right.\}$ by 
\[
\begin{pmatrix} a & b\\c & d\end{pmatrix} z = \dfrac{a z + b}{c z + d}.
\]
Recall also that we can define a metric on $H$ by setting 
\[
d(z, w) = \acosh\left( 1 + \frac{|z-w|^2}{2 \Im z \Im w}\right),
\]
and, equipped with this metric, $H$ is isometric to the hyperbolic plane $\mathbb{H}^2.$ In addition, the action of $\SL(2, \reals)$ by linear fractional transformations described above is isometric, and, indeed, the every isometry of $\mathbb{H}^2$ is obtained this way, so 
\[
\Isom \mathbb{H}^2 \simeq P\SL(2, \reals) = \SL(2, \reals)/\{\pm I\},
\]
where the quotient by plus and minus identity is needed because $(-I) z = \frac{-z}{-1} = z,$ for all $z \in H.$
Weewill also need the \emph{singular value decomposition}. Recall that every matrix $A$ in $M^{m \times n}$ can be written as $A=P D Q,$ where $P \in O(m),$ $Q \in O(n),$ and $D$ is diagonal $m\times n$ matrix with nonnegative diagonal elements (see, e.g., \cite{horn1990matrix}). The diagonal elements of $D$ are known as the $\emph{singular values}$ of $A.$ It is well-known (and easy to verify) that the Frobenius norm of $A$ equals the Euclidean ($L^2$) norm of the vector of its singular values.

In the special case where $n=m=2,$ and $\det A = 1,$ it is easy to see that the above implies that $A$ can be written as 
\[A=\begin{pmatrix} \cos \phi & \sin \phi\\ -\sin \phi & \cos \phi \end{pmatrix} \begin{pmatrix} x & 0 \\ 0 & \frac1x\end{pmatrix}\begin{pmatrix} \cos \theta & \sin \theta\\ -\sin \theta & \cos \theta \end{pmatrix}, \]
for some $x > 1.$ Further, as noted above, $\|A\|^2 = x^2 + 1/x^2.$
\subsection{Translation distance} A big part of the reason for introducing the singular value decomposition above is to give a palatable answer to the following question:
\begin{question}
\label{distq}
How far (in hyperbolic metric) does the matrix $A = \begin{pmatrix}a & b \\ c & d\end{pmatrix} \in \SL(2, \reals)$ move the point $i?$
\end{question}
The main reason why the singular value decomposition helps is that 
\[\begin{pmatrix} \cos \theta & \sin \theta\\ -\sin \theta & \cos \theta \end{pmatrix} i = i,\]
so with $A$ as above, we have 
\[
A i = \dfrac{-x^2 \cos \phi + i \sin \phi}{i \cos \phi + x^2 sin \phi}.
\]
After some tedious computation (or a couple of lines of \textit{Mathematica}) we obtain:
\begin{gather}
\Re Ai = \dfrac{\cos \phi \sin \phi}{\cos^2 \phi + x^4 \sin^2 \phi} - \dfrac{x^4 \cos\phi \sin\phi}{\cos^2 \phi -+ x^4 \sin^2 \phi} \\
\Im Ai = \dfrac{x^2}{\cos^2 \phi + x^4 \sin^2 \phi},
\end{gather}
and finally
\begin{equation}
\label{disteq}
d(i, Ai) = \acosh\left(1+\dfrac{(x^2-1)^2}{2x^2}\right) = \acosh\left(\frac12\left[\frac1{x^2} + x^2\right]\right) = 2\log x,
\end{equation}
a surprisingly simple answer, after all that computation.

As a minor bonus, we can now modify our procedure $\PickHyperbolic$ to return a point in the upper halfplane in procedure $\PickHalfplane$ (see Algorithm \ref{pickhalfplane}).
\begin{algorithm}
\label{pickhalfplane}
\caption{Picking a random point in the disk around $i$ in the Poincar\'e halfplane model}
\begin{algorithmic}[1]
\Require $R$ a positive real number.
\Function{PickHalfplane}{$R$}
\Statex
\State $(r, \theta) \gets \PickHyperbolic(R)$
\State \Return $\dfrac{i \cos e^x \theta  \sin \theta}{-i e^x sin \theta + \cos \theta}$
\EndFunction
\end{algorithmic}
\end{algorithm}

\subsection{The fundamental domain and orbits of the $\SL(2, \integers)$ action}
\begin{figure}[hb]
\label{moddom}
\centering
\includegraphics[width=4in]{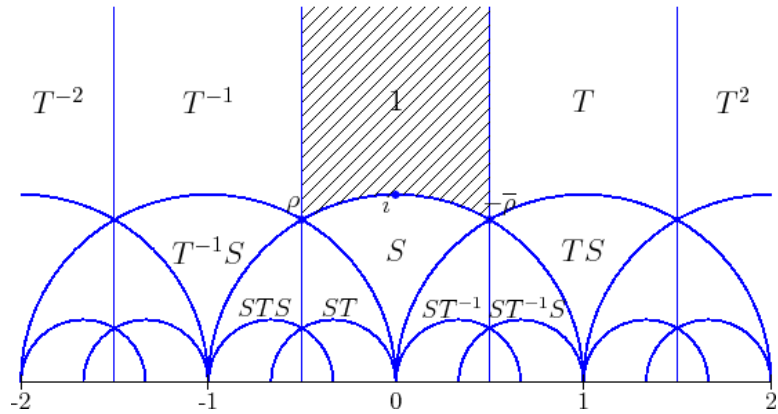}
\caption{The modular tesselation; fundamental domain shaded, other copies labeled by the elements sending the shaded domain to the copy}
\end{figure}

The action of $\SL(2, \integers)$ on $H$ is discrete, and its fundamental domain $\Lambda$ is one of the best known images in all of mathematics (the reader can see it again in Figure \ref{moddom}). 
The points in the fundamental domain index the orbit of the $\SL(2, \integers)$ action, and gives rise to the following natural question:
\begin{question}
\label{whichorbit}
Given a point $z\in H,$ which orbit is it in? In other words, which point of $\Lambda$ gets mapped to $z?$
\end{question}
This question is so natural it was asked and answered in the 18th century by Legendre and Gauss. Of course, for them, the question was a little different: they were given two linearly independent vectors in the plane. These vectors generate a lattice, and the question is: what is the canonical form for that lattice? In other words, Gauss and Legendre posed (and solved) the two dimensional \emph{lattice reduction} problem (a very nice reference is the paper \cite{vallee2007lattice}). Gauss' algorithm (which is basically the continued fraction algorithm) proceeds as follows:
\begin{algorithm}
\label{latticered}
\begin{algorithmic}[1]
\Require A complex number $z$ with $\Im z \geq 0.$
\Function{Reduce}{$z$}
\While{$z\leq1$}
\State $z \gets -1/z$
\State $q \gets \round{\Re z}$
\State $z \gets z - q$
\EndWhile
\State \Return $z$
\EndFunction
\end{algorithmic}
\end{algorithm}
In fact, Algorithm Reduce can be made to do more: give the point $z\in H,$ we can return not just the point $z_0 \in \Lambda$ such that $z$ is in the orbit of $z_0,$ but also the matrix $A \in \SL(2, \integers)$ such that $z_0 = A z,$ as done in Algorithm Reduce2.
\begin{algorithm}
\label{latticered2}
\begin{algorithmic}
\Require A complex number $z$ with $\Im z \geq 0.$
\Function{Reduce2}{$z$}
\State $A \gets I$
\While{$z\leq1$}
\State $z \gets -1/z$
\State $A \gets \begin{pmatrix}-1 & 0 \\ 0 & 1\end{pmatrix} A$
\State $q \gets \round{\Re z}$
\State $z \gets z - q$
\State $A \gets \begin{pmatrix}1 & q \\ 0 & 1\end{pmatrix} A$
\EndWhile
\State \Return $(A, z)$
\EndFunction
\end{algorithmic}
\end{algorithm}
\section{Selecting a random element of $\SL(2, \integers)$ almost uniformly.}
We are now ready to describe the algorithm for selecting a random matrix $M$ from the set  of matrices in $\SL(2, \integers)$ with Frobenius norm bounded above by $X.$ Aside from the observations above, the key remark is that the Haar measure on $\SL(2, \reals)$ projects to the hyperbolic metric on $H$ (see the discussion in \cite{eskin1993mixing,duke1993density}). This suggests the following algorithm:
\begin{algorithm}
\label{randmatsl2}
\caption{Returns a matrix in $\SL(2, \integers)$ with Frobenius norm bounded by $X.$ The ratio of the probabilities of any two matrices s between $e^\epsilon$ and $e^{-\epsilon}$}
\begin{algorithmic}[1]
\Require A pair of positive real number $X,\epsilon$
\Function{PickFancy}{$X$,$\epsilon$}
\State $R\gets f(X,\epsilon)$ \Comment $f$ is a function to be named later
\Loop
\State $z \gets \PickHalfplane(f(X))$
\State $(A, z_0) \gets \Reduce2(z)$
\If{$\|A\|\leq X$}
\State \Return $A$
\EndIf
\EndLoop
\EndFunction
\end{algorithmic}
\end{algorithm}

What should $f(X, \epsilon)$ be? Firstly, it is obviously necessary that the disk of radius $f(X, \epsilon)$ intersect \emph{all} of the fundamental domains of matrices $A$ As we have seen (Eq. \eqref{disteq}), in order for this to be true, we must have $f() > 2 \acosh 2 X.$ On the other hand, the fundamental domain $\Lambda$ of $\SL(2, \integers)$ has a cusp, which is bad, since no disk can contain $\Lambda,$ but not \emph{so} bad, since the part of $\Lambda$ which lies outside the disk of radius $R$ around $i$ is asymptotic to $ \exp(-R+1).$  This means that if $f(X) > t+ \acosh 2 X^2,$ the ratio of the areas of the intersections of fundamental domains we are interested in is of order $1+e^-t.$ On the other hand, the number of fundamental domains we \emph{do not} want is proportional to $e^t,$ so, as claimed in the introduction, the amount of excess computation is proportional to the error.
\subsection{Complexity estimates and implementation}
\label{complexity}
Picking the random number in the halfplane in function $\PickHalfplane$ has been made unnecessarily expensive. Unwinding what we are doing, we see that in the first step we 
pick  a random number $x$ between $0$ and $\cosh(C+\acosh 2 X^2) - 1,$ which is an algebraic function of $X,$ and in the next step we generate the $\acosh(x + 1)$, which is  a combination of logarithm and square root. Since the number of fundamental domains is exponential in the radius, we need roughly $\log X$ bits of precision, and the final step ($\Reduce2$) then takes a logarithmic number of steps (see \cite{lagarias1980worst}), each of which is of logarithmic complexity (note that Daub\'e, Flajolet, and Valle\'e \cite{daude1996average} show that with the uniform distribution, the expected number of steps does \emph{not} depend on the size of the input, but it remains to be investigated whether this is true for our model).
\section{Extensions to other Fuchsian and Kleinian groups}
Suppose that instead of $\SL(2, \integers)$ we want to generate random elements of bounded norm from other subgroups of $\SL(2, \reals)$ or, even more ambitiously, $\SL(2, \cx).$ The general approach described above works. Suppose $H$ is our (discrete) subgroup. To pick a random element, we pick a random point $x$ in $\hyps^2$ or $\hyps^3,$ (our radius computation goes through unchanged) then find the matrix $A \in H$ which moves $x$ to the ``canonical'' fundamental domain of $H.$ This last part, however, is not so obvious, because both questions (constructing the fundamental domain and ``reducing'' the point $x$ to that fundamental domain) are nontrivial.

\subsection{Constructing the fundamental domain} The first observation is that if the group $H$ is not geometrically finite, it does not have a finite-sided fundamental domain at all, so constructing one may be too much. It is, however, conceivable that deciding whether $x$ is reduced (that is, lies in the canonical fundamental domain) is still decidable. Since no algorithm leaps to mind, we shall state this as a question:
\begin{question}
\label{geomfin}
Is there a decision procedure to determine whether $x \in \hyps^n$ lies in the canonical fundamental domain for a not-necessarily-geometrically finite group $H?$
\end{question}
Until Question \ref{geomfin} is resolved, we will assume that $H$ is geometrically finite. Now, we can construct the fundamental doman by generating a chunk of the orbit of the basepoint, and then computing the Voronoi diagram of that pointset -- the resulting domains are the so-called \emph{Dirichlet} fundamental domains. Computing the Voronoi diagram can be reduced to a Euclidean computation (see the elegant exposition in \cite{nielsen2010hyperbolic}, and H. Edelsbrunner's recent classic \cite{edelsbrunner2001geometry} for background on the various diagrams). However, a much harder problem is of figuring out how much of an orbit needs to be computed. For Fuchsian groups, this was addressed by Jane Gilman in her monograph \cite{gilman1995two} (at least for two-generator fuchsian groups). For Kleinian groups the question is that much harder, but has been studied at least for \emph{arithmetic} Kleinian groups in \cite{page2012computing}. All we can say in general is that the computation is finite (since at every step we check the conditions for the Poincar\'e polyhedron theorem), so after waiting for a finite (though possibly long) time, we are good to go. Now, the question is: lacking the number theory underlying the continued fraction algorithm, how do we \emph{reduce} our random point to the canonical fundamental domain? There are a number of ways to try emulate the continued fraction algorithm. Here is one.
\begin{algorithm}
\label{tryreduce1}
\caption{greedy reduction algorithm}
\begin{algorithmic}
\Require $x,b \in \hyps^n,$ side-pairing transformation of the Dirichlet domain $\Gamma=\{ \gamma_0=I(n),\gamma_1, \dotsc, \gamma_k\}$
\Function{GreedyReduce}{$x$, $b, \Gamma$} \Comment $b$ is the basepoint.
\State $M \gets I(n)$
\Loop
\State Loop over $\Gamma$ to find the $i \in [0, k]$ for which $d(\gamma_i(x), b)$ is minimal.
\If{$i=0$}
\State \Return $M$
\EndIf
\State $M \gets \gamma_i M$
\State $ b \gets \gamma_i b$
\EndLoop
\EndFunction
\end{algorithmic}
\end{algorithm}

Algorithm \ref{tryreduce1} will terminate in \emph{at most} exponential time (exponential in $d(b, x),$ that is), and it seems very plausible (for reasons of hyperbolicity) that it will actually terminate in time \emph{linear} in $d(b, x),$ but this seems difficult to show.

\section{Higher rank}
\subsection{$\SL(n, \integers)$}
The algorithms for $\SL(2, \integers)$ use, in essence, the $KAK$ decomposition of the group (which is in this case the singular value decomposition). This exists, and is easy to describe geometrically, in the higher rank case as well (this construction is due to Minkowski). We first introduce the \emph{positive definite cone} \[\PSD(n) = \{M \left| M=M^t, v^t M v \geq 0, \forall v \in \reals^n\right.\}\] The general linear group $\GL(n,\reals)$  acts on $\PSD(n)$ by $g(M) = gMg^t$ It is not immediate that the subset $\PSD_1(n) = \{M \in \PSD(n) \left| \det M = 1\right.\}$ is invariant under $\SL(n, \reals)$ We can define a family of (Finsler) metrics on $\PSD(n)$ by \[d_p(A, B) = \left(\sum_{i=1}^n| \log \sigma_i(B^-1 A)|^p\right)^{1/p},.\] where $\sigma_i(M)$ denotes the $i$-th singular value of $M.$ When $p=2$ this defines a Riemannian metric, which makes $\PSD_1(n)$ into the symmetric space for $\SL(n, \integers).$ In particular, when $n=2$ it is easy to check that $\PSD_1(2)$ the hyperbolic plane $\hyps^2$ with the usual metric. With this in place, the algorithm we described for $\SL(2, \integers)$ goes through \emph{mutatis mutandis}. The hard part is the reduction algorithm. In the setting of $\SL(n, \integers)$ we have the \emph{lattice reduction} problem, which has been heavily studied starting with L. Lovasz' foundationalLLL algorithm in \cite{lenstra1982factoring}. The LLL algorithm is generally used as an \emph{approximation} algorithm: it reduces a point not into the fundamental domain but into a point near the fundamental domain, which begs the question:
\begin{question}
\label{llldist}
Are the matrices obtained in the LLL algorithm uniformly distributed?
\end{question}
In any case, one can also perform \emph{exact} lattice reduction, but in that case the running time is exponential in dimension (sse \cite{nguyen2011lattice}); for dimensions up to four there is an extension of the Legendre-Gauss algorithm, described above, which is exact and quadratic in terms of the bit-complexity of the input, see \cite{nguyen2009low}.
\subsection{$\Sp(2n, \integers)$}
For $\Sp(2n, \reals)$ the symmetric space is the Siegel half-space, where the metric is defined the same way as for $\SL(n, \reals),$ while the underlying space is not the positive semidefinite cone, but instead the set $S(2n)$ f all \emph{complex} symmetric matrices with poisitive definite imaginary part. A symplectic matix $X \in \Sp(2n, \reals)$ has the form $X=\begin{pmatrix}A & B\\ C & D\end{pmatrix}$ where $A, B, C, D$ are $n\times n$ matrices satisfying the conditions that $A^t C A^{-1} (C^{-1})^t = B^t D B^{-1} (D^{-1})^t = A^t D - C^t B = I(n).$ The action of $\Sp(2n, \reals)$ on $S(2n)$ is then given by:
\[ X(Z) = (AZ + B)(CZ + D)^{-1}.\]
For more details on this, see \cite{siegel1943symplectic,freitas1999action}. In any case, the action of $\Sp(2n, \integers)$ on the Siegel half-space is fairly well understood, and the algorithm we gave for $\SL(2, \integers)$ (which is also known as $\Sp(2, \integers)$) goes through, with the usual question of lattice reduction, which has \emph{not} been studied very extensively; the only reference I have found was \cite{gama2006symplectic}, which is, however, quite throrough.

\section{Miscellaneous other groups}
\subsection{The orthogonal group}
Even without integrality assumptions, it is not immediately obvious how to sample a uniformly random matrix from the orthogonal group. This question got a very elegant one-line answer from G. W. Stewart in his paper \cite{stewart1980efficient}. Stewart's basic method is as follows: Firstly, we remark that it is well-known that every matrix $M$ possesses a $QR$ decompoosition, where $Q$ is orthogonal, while $M$ is upper triangular, and this decomposition is unique up to post-multiplying $Q$ by a diagonal matrix whose elements are $\pm 1.$ This indeterminancy can be normalized away by requiring the diagonal elements of $R$ to be positive. The algorithm is now the following(Algorithm \ref{stewartalg}):
\begin{algorithm}
\label{stewartalg}
\caption{Generating matrices in $\O(n)$ uniform with respect to the Haar measure}
\begin{algorithmic}
\Require $n$ is a positive integer.
\Function{RandomOrthogonal}{$n$}
\State $X \gets \mbox{an $n\times n$ matrix whose entries are independent with the common distribution $N(0, 1)$}$
\State $(Q, R) \gets \mbox{the QR decomposition of $X$}$
\State \Return $Q$
\EndFunction
\end{algorithmic}
\end{algorithm}
This algorithm works because the distribution of $K X$ is the same as the distribution of $X$ for a matrix $X$ with i.i.d. normal entries, and so the distribution of $KQ$ is the same as the distribution of $Q,$ which is exactly what we seek (notice that this method is morally a slight extension of the method described in Section \ref{randomballs}), and is also morally related to our algorithms for $\SL(n, \integers).$

 Now generating random \emph{integral} matrices in $O(n)$ is easy -- they are just the signed permutation matrices, and generating a random permutation is easy (in a quest for self-containment we give the algorithm below as Algorithm \cite{permalg}, as is assigning random signs. However, as far as I know there is no known way to generate uniformly random \emph{rational} orthogonal matrices. We ask this as a question:
\begin{question}
\label{orthog}
How do we generate a random element of $O(n)$ whose elements have greatest common denominator bounded above by $N?$
\end{question}

There is a natural companion question:
\begin{question}
\label{orthog2}
Let $O_q(n)$ be the set of those elements of $O(n)$ with rational entries, such that the size of the greatest common denominator is bounded above by $q.$ Is there any exact or asymptotic formula for the order of $|Q_q(n)|?$
\end{question}

And another natural question:
\begin{question}
\label{orthog3}
Let $\mu_q$ be the normalized counting measure on $O_q(n)$ (as above). Do the measures $\mu_q$ converge weakly to the Haar measure on the orthogonal group?
\end{question}

Questions related to Questions \ref{orthog2} and \ref{orthog3} are considered in the paper \cite{MR2482443}, and it is quite plausible that the methods extend,  but it is not completely obvious as of this writing. The only thing we know with certainty is how to address the case of $SO(2).$ Here, the elements have the form $\begin{pmatrix} a & b \\ -b & a\end{pmatrix},$ with $a^2 + b^2 = 1.$ Thus, if $a$ and $b$ have denominator $q,$ we are counting the representations of $q$ as a sum of two squares. For this there is the explicit formula of Dirichlet:
\begin{quotation}
If $q = p_1^{2a_1} \dots p_k^{2a_k} q_1^{b_1} \dots q_l^{b_l},$ where $p_i = 4k_i + 3,$ which $q_j= 4k_j + 1,$ then the number of way to write $q$ as a sum of two squares is $\prod_{j=1}^l (b_j + 1).$
\end{quotation}
To get an asymptotic result, it is necessary to consider all $q \leq Q,$ when we see that the number of elements in $SO(2)$ with the greatest common divisor of coefficients equals the number of \emph{visible} lattice points in the disk $\|x\| \leq Q$ (a visible point $(a,  b)$ is a lattice point with relatively prime $a, b$). Since the probability of a lattice point being relatively prime for $Q \gg 1$ approaches $6/\pi^2,$ and the number of lattice points in the disk is asymptotic to $\pi Q^2,$ we see that the cardinality of $SO_Q(2)$ is asymptotic to $\frac{6}{\pi} Q^2,$ so we have a rather satisfactory answer to Question \ref{orthog2} in this setting. 

Question \ref{orthog3} is also easy (but already deep) in this setting. It is equivalent to the equidistribution of rational numbers with bounded denominator in the interval, and that, it turn, is not hard to show is equivalent to the prime number theorem (both statements are equivalent to the statement that $\sum_{k=1}^x \mu(x) = o(x),$ where $\mu$ is the M\"obius function).

Finally, in view of the answer to Question \ref{orthog2}, Question \ref{orthog} is equivalent to the question of generating a lattice point in a ball, which we have already discussed in Section \ref{randomballs}
\begin{algorithm}
\label{permalg}
\caption{generating a random permutation uniformly}
\begin{algorithmic}[1]
\Require $n > 0$
\Function{GenPerm}{$n$}
\State $a \gets [1, 2, 3, \dots, n]$
\For{$i=1 \to n$}
\State \mbox{swap $a[1]$ and $a[n-i+1]$}
\EndFor.
\State \Return $a$
\EndFunction
\end{algorithmic}
\end{algorithm}
\subsection{Finite Linear Groups}
Our final remarks are on finite linear groups. The simplest class of groups to deal with is $\SL(n, p)$ How do we get a random element? This is quite easy, see Algorithm \ref{simpleslnfinitealg}:
\begin{algorithm}
\label{simpleslnfinitealg}
\caption{generating a random element of $\SL(n, p).$}
\begin{algorithmic}[1]
\Require $n > 0$
\Function{GenRandSL}{$n$}
\Loop
\State $a \gets \mbox{a uniformly random element of $M^{n\times n}(p).$}$
\If{$\det(a) \neq 0$}
\State
\Return $a$ with the first column divided by $\det(a).$
\EndIf
\EndLoop
\EndFunction
\end{algorithmic}
\end{algorithm}
We  pick every element independently at random from $F_p.$ If the resulting matrix $M$ is singular, we try again, if not, let the determinant  be $d.$  We then divide the first column of $M$ by $d.$ It is easy to see that the resulting matrix $M^\prime$ will  be uniformly distributed in $\SL(n, p).$ It is easy to see that the complexity of this method is $O(n^\omega \log p),$ where $\omega$ is the optimal matrix multiplication exponent.Unfortunately, this simple method only works for $\SL(n, q).$ For $\Sp(2n, q)$ there is the Algorithm \ref{hallalg}, which is due to Chris Hall.
\begin{algorithm}
\label{hallalg}
\caption{Chris Hall's algorithm to generate a random element of $\Sp(2n, p).$}
\begin{algorithmic}[1]
\Require $n > 0$
\State {$V \gets \mbox{symplectic vector space of dimension $2n.$}$}
\Function{GenRandSp}{$n$}
\State{$W\gets \{0\}$}

\For{$i=1 \to n; i \gets i+1$}
\Repeat
\State {$x, y \gets \mbox{random vectors in $V.$}$}
\State{$x^\prime, y^\prime \gets \mbox{projections of $x, y$ onto $W.$}$}
\State{$x^{\prime\prime} \gets x - x^\prime$}
\State{$c\gets \langle x^{\prime\prime}, y^{\prime\prime}\rangle$}
\Until{$c \neq 0$}
\State{$x_i\gets x^{\prime\prime}$}
\State{$y_i\gets y^{\prime\prime}/c$}
\State{$W\gets \mbox{span of $W$ and $x_i, y_i$}$}
\EndFor
\State
\Return{$x_1, x_2, \dotsc, x_n, y_1, y_2, \dotsc, y_n$}
\EndFunction
\end{algorithmic}
\end{algorithm}
It is not hard to see that Chris Hall's algorithm has time complexity $O(n^3 \log p).$

In general, there is a completely different polynomial-time algorithm based on the fact that the Cayley graphs of simple groups of Lie type are expanders -- uniform expansion bounds have been obtained by a number of people, see \cite{MR2746060,MR2342638,MR2342639,lubexpanders} The main significance of the expansion for our purposes is tha the random walk on the Cayley graph is very rapidly mixing -- see \cite[Section 3]{hoory2006expander}, and so a random walk of polylogarithmic length will be equidistributed over the group. Of course, this will be slower than Algorithm \ref{simpleslnfinitealg} , and will only generate \emph{approximately} uniform random elements. To be precise, the diameter of the Cayley graph of (for example) $\SL(n, p)$ will be $O(n^2 \log p),$ so the expander-based algorithm will have time complexity $O(\log^2 p n^{\omega + 2}).$

\subsection{Other groups?}
In the work by the author \cite{rivin2008walks,rivin2009walks} and Joseph Maher (\cite{MR2772067}) the model of a random element is the random walk model, since this seemed to the only natural model for the mapping class group. However, in view of the discussion above it makes sense to define the norm of an element $\gamma$ of a mapping class group as the Teichmuller distance from some fixed base surface $S$ to $\gamma(S)$ (one can also use the Weil-Petersson distance, or the distance from a fixed curve to its image in the curve complex, and then pick a random element by analogy with the construction in this note. In fact, this has been done by Joseph Maher in \cite{maher2010asymptotics}. 

\bibliographystyle{plain}
\bibliography{msri}
\end{document}